\documentclass[a4paper,UKenglish]{lipics}
 
\usepackage{microtype}

\title{Making Room for Women in our Tools for Teaching Logic: A Proposal for Promoting Gender-Inclusiveness\footnote{This work was partially supported by a Capes Postdoctoral Researcher Grant.}}
\titlerunning{Making Room for Women in our Tools for Teaching Logic} 

\author[1]{Frederique Janssen-Lauret}
\affil[1]{Department of Philosophy, University of Campinas\\
   Rua Cora Coralina 100, CEP 13083-896 Campinas, Sao Paulo, Brazil\\
  \texttt{fmjanssenlauret@cle.unicamp.br}}

\authorrunning{F.\,M. Janssen-Lauret} 

\Copyright{Frederique Janssen-Lauret}


\serieslogo{logo_ttl}
\volumeinfo
  {M. Antonia {Huertas}, Jo\~ao {Marcos}, Mar\'ia {Manzano}, Sophie {Pinchinat}, \\
  Fran\c{c}ois {Schwarzentruber}}
  {5}
  {4th International Conference on Tools for Teaching Logic}
  {1}
  {1}
  {65}
\EventShortName{TTL2015}

\begin{document}

\maketitle

\begin{abstract}
Logic is one of the most male-dominated areas within the already hugely male-dominated subject of philosophy. Popular hypotheses for this disparity include a preponderance of confident, mathematically-minded male students in the classroom, the historical association between logic and maleness, and the lack of female role-models for students, though to date none of these have been empirically tested. In this paper I discuss the effects of various attempts to address these potential causes whilst teaching second-year formal and philosophical logic courses at different universities in the UK. I found the most noticeable positive effect came from assigning a good proportion of reading by female authors presenting an original point of view. I go on to suggest some implementations for incorporating more texts by female authors into our arsenal of tools for teaching logic.
 \end{abstract}

\section{Introduction}

Logic is a heavily male-dominated discipline. Men prevail on reading lists, at conferences, and in journals. Although this fact cannot have escaped notice, even those of us who are aware of it and actively attempt to counteract it may be temporarily stunned when presented with the evidence of just how rare it is to see due attention being paid to a woman's original argument in a teaching context in formal or philosophical logic. 

Consider the following illustrative example. The University of Cambridge second-year Logic reading list---which covers both formal and philosophical logic---contains, at the time of writing, work by ninety-four distinct men, some of whom are cited many times over; by contrast, only seven women are listed \cite{cambridgereadinglist14}. Two of these women have male co-authors. Four of them are the sole authors of textbooks, historical books, or survey papers; they are presented primarily as commentators on the work of men.  Only one of the papers recommended to second-year logic students at Cambridge is a piece of original research in which a female author presents her own point of view: Dorothy Edgington on the paradox of knowability \cite{edgington85}.

Men dominate not only the literature,  but also the classroom, both as teachers and students.  As lecturers and tutors on courses in formal and philosophical logic, we know that these courses generally contain more male students than female students. This appears to be partly due to the general gender disparity in subjects like philosophy and mathematics, whose students generally take logic courses, but is also seen to be a result of self-selection by female students away from the area of logic and towards more normative and applied topics. We also frequently see that the male students in courses on formal and philosophical logic are much more vocal than the women. And even though philosophy is generally male-dominated \cite{beebeesaul11}, this  disparity between male and female students appears to be greater in the area of logic than in other areas of philosophy. Unfortunately no systematic research has been carried out into the reasons for this gender imbalance. But in discussions on the subject, we often hear concerned people point to one or more of the following: the dominance of confident, mathematically-minded male students in the classroom, the historical association between logic and maleness, and the lack of female role-models for students. 

In this paper I discuss the results of various strategies which I have used to address these possible causes of gender inequality while teaching formal and philosophical logic to second-year students at the University of St Andrews and at four Cambridge colleges: Clare, Newnham, Selwyn, and Trinity Hall. My ability to make significant changes and investigate issues systematically within these courses was quite limited, because I was involved only in a teaching assistant capacity. But I set out to counteract the potential sources of inequality as best I could with the resources I had available to me. To dispel the impression that men have more of an aptitude for logic or mathematics, I emphasised the importance of working hard over innate ability, encouraged students' efforts in formal work without praising mathematical talent, and drew their attention to the continuity between formal logic and philosophy. All of these, I found, anecdotally had some positive effect in terms of either student feedback or test scores, but the most dramatic positive effect came from consistently assigning a good proportion of original research by female authors as part of the set reading. My proposal is therefore to build teaching materials, including syllabuses, reading lists, and textbooks more specifically around the ideas of women. Where possible a female author should always be presented as the originator of an important view of her own, rather than as a survey author commenting on the work of men.

\section{Counteracting Chauvinism and the Association Between Logic and Maleness}

When I talk about a pervasive association between logic and maleness as a potential reason that women are put off the academic study of logic, I am not talking about anything as all-encompassing as the idea that the methods of logic are intrinsically linked with oppressive discourse, as proposed, for instance, by Andrea Nye \cite{nye90}. I mean the idea, or stereotype, which is still commonly held either implicitly or explicitly, that men naturally have more of an aptitude for formal logic and mathematics, or even for rigorous thinking and rationality in general.\footnote{Women are not the only marginalised group subjected to the harmful perception that they are less rational than some other, socially dominant group of people. There are racist versions of this stereotype as well, and some currently popular arguments for atheism also paint all religious people with the irrationality brush. Some students, especially Muslim students and others belonging to minority religions, have told me that this has led to their feeling alienated from mainstream philosophy. I am not aware of systematic research on this specific question, but I hope that some of the strategies outlined in this paper may also be effective against these kinds of harmful stereotypes. (As I make clear below, because of the risk of stereotype threat I deliberately avoided saying explicitly that these were efforts to promote gender-inclusiveness specifically.)} It is likely that such stereotypes have some effect on female students that results in turning them away from courses in logic, and cause them to be more quiet and diffident compared to their male counterparts when they do take such courses. While teaching formal logic exercise classes in an intermediate course at the University of St Andrews---comprising some model theory, natural deduction, and Kripke models for modal logics up to S5---I made efforts to challenge this idea. Anecdotally this led to some overall improvement in exam results, and it also occasioned some positive comments on student questionnaires. 

\subsection{No Good Evidence That Men Have Greater Mathematical Ability} First of all it is helpful to consider the figures which show that the stereotype really is false. There is no clear evidence to corroborate the idea that men have more of an aptitude for mathematics than women do, nor are there recent data that clearly demonstrate that men perform significantly better at mathematical tests than women do. Although in some countries male pupils significantly outperform female pupils, the data available world-wide, for instance those from the Programme for International Student Assessment (PISA), yield very mixed results. Several other countries have female pupils outperforming their male counterparts, or have no gender difference \cite[pp. 71-75]{pisa12}. A large global meta-analysis also found only very minor overall gender differences, but great variation between different countries \cite{elsequesthydelinn10}. The causes of this variation are unknown, and there is no clear hypothesis; high performance in mathematics for female pupils does not correlate especially strongly with higher status for women generally \cite[p. 411]{bank07}. Different analyses also reveal different results: although the US is one of the countries where according to PISA, male pupils outperform female pupils on tests, large US-wide meta-analyses have shown no significant differences in school performance in mathematics between male and female pupils \cite{hyde08, hyde10}. Overall the results appear consistent with the gender similarities hypothesis:  that the psychology of men and women is much more similar than different \cite{hyde05}, equally so with respect to mathematical ability.

\subsection{Avoiding Stereotype Threat} Nevertheless, we know that stereotypes, even when false, are themselves causally efficacious in negatively affecting the test performance of people in underrepresented groups. When teaching intermediate formal logic to second-year students, I wanted to make explicit efforts to counteract the known effects of stereotype threat: the phenomenon of underperformance in disadvantaged groups when they are reminded of their disadvantaged status. This has been shown to apply to women's maths performance \cite{spencer99, mcglone07}. Since formal logic courses are perceived by students as being very mathematical in nature, and are sometimes taken by mathematics students alongside philosophy students (as was the case in the formal logic course I taught), this kind of stereotype threat is likely to affect female students taking courses on formal logic. This implies, amongst other things, that well-intentioned efforts to assure female students explicitly that they are just as capable as the male students have the potential to backfire. 

\subsection{Emphasising Motivation and Effort over Innate Talent}  Another issue that has been identified by social and developmental psychologists as having the potential to hold back female students is that of gendered expectations concerning motivation and feedback. There is some evidence, deriving from extensive studies by Carol Dweck amongst others, that praising students for their efforts, rather than praising them for their innate ability, generally produces more improvement in their performance. But the overall trend is complicated by gendered teaching environments. In Dweck's studies from the 1970s and 80s, teachers were found to attribute male pupils' failures more often to lack of motivation or situational factors, and girls' failures to intellectual inadequacies \cite{dweck78}. She hypothesised that this might account for gender differences in mathematics achievement, which were much more pronounced in studies from this period than they are now \cite{dweck86}. 

While the gender gap in school pupils' mathematics achievement has now narrowed, this effect may well still be at work in university-level logic teaching. More recent work by Dweck indicates that adult women are still more likely, in the twenty-first century, to believe mathematics to be a matter of innate ability, and to ascribe their own lack of mathematical achievement to intellectual inadequacy \cite{dweck06}.  There are no studies in this tradition that pertain specifically to undergraduate logic courses, but my own experience suggests that the stereotype is persistent. I have certainly personally heard people, including some of my students, say that men are simply better at logic. In my logic classes, I've frequently witnessed  female students (and some male students, too) give up at an early stage complaining they don't have `that kind of brain', attributing their difficulties to lack of innate ability. Male students' difficulties are chalked up to lack of trying, dyslexia, \emph{etcetera}, their successes to mathematical ability. Female students invested in their self-perception as having `the wrong kind of brain' will explicitly attribute other women's good performance in formal logic to their being atypical women. It may seem as though having female tutors and lecturers would dispel this perception, but as a female logic tutor, and now a female academic working on philosophical logic, I have not found this to be consistently true. I have been told many times, sometimes by female as well as male students, that I must just not be a typical woman to be good at or enjoy logic.\footnote{I am not alone in this; women who are good at any area of the stereotypically masculine subject of philosophy are routinely told that this means they must be masculine somehow \cite[p. 212]{haslanger08}.} I have never had an explicit positive comment from a student about having a female tutor. Specific strategies are needed to make logic-teaching more gender-inclusive.

\subsection{My Strategies for Addressing Male-Dominated Classrooms and Harmful Stereotypes} To dispel the impression that logic is a subject that men are naturally more suited to than women, I tried to apply Dweck's theory by emphasising the importance of working hard over innate ability. I wanted to do this in a way that did not explicitly reference the stereotype that logic is `for men', since this itself carries with it the risk of stereotype threat, which is observed whenever people in disadvantaged groups are explicitly reminded of their disadvantaged status. Instead I attempted to make positive efforts to include all students in completing, demonstrating, and explaining the exercises, and to make them all feel like the course was one they could do well at with some effort, not one that divided the group into the talented ones who would excel at it and the untalented ones who had no hope of doing so.

I stressed frequently to the students that success in logic at the undergraduate level is \emph{not} a matter of talent, or `having a head for it', or `having that kind of brain', to use the phrases students often repeat to each other which I view as reinforcing their exaggerated impression of the role played by innate ability. I explained that they would find that the secret to doing well in a logic course is to keep on top of the material and do the exercises every time, without fail. I emphasised that most students were perfectly capable of handling the material, and indicated where to find extra help for those who found it challenging, including help for specific learning difficulties like dyslexia or dyscalculia. I wanted to move their thinking away from the rather disablist idea that specific learning difficulties should be classed as `having the wrong kind of brain', when universities are perfectly capable of providing support and accommodations for these issues. Some students with learning difficulties were provided with accommodations such as extra exam time and dictation software, and did well in the course as a result. 

I also explicitly gave the students strategies for reviewing material covered earlier in the course if they were lost or felt baffled by an exercise. In particular, tracing back their steps to last week's exercises, identifying the new component, regaining their confidence in their ability to handle last week's material, and tackling the difficult question with a specific understanding of its relationship to familiar issues and an ability to isolate the newer or more challenging components worked very well for many of the more logic-averse students. 

Of course, these solutions are effective only where students are actually made to do the work consistently.  I made it mandatory for all students to hand in a full set of proofs every week. This required them to submit their best attempt at solving the problem. I did not allow them to skip any questions or claim that some of the exercises were beyond them. They had to submit at least a few lines of each proof before they could be counted as present for the class, and be given points for attendance. After the initial mild shock all students complied with this policy. I then asked all students in turn to present their proofs on the whiteboard, explain their reasoning to each other, and help each other find the correct answer, stepping in from time to time to moderate and ensure that the discussion was not dominated by the more confident students. 

I also made efforts to point out possible applications of the principle at work in a proof, or continuities between logic and philosophy. Some of the more logic-averse students told me that this made them view the course as more interesting and relevant to them. 

\subsection{Summary of Positive Effects} Although I did not get the chance to investigate the effects of my strategies here in a systematic way, they did appear to lead to improved confidence in previously logic-averse students, some of them male, some female.\footnote{I have since become aware of some promising empirical work that is being done on inclusive methodology in philosophy teaching, though I haven't found anything specific to formal or philosophical logic. For instance, Kristina Gehrman is currently conducting a multi-year study in her introductory philosophy course at the University of Tennessee \cite{gehrman15}; her  methodology is based on social-role congruence models familiar from science and mathematics \cite{diekman13}. And Eva Cadavid at Centre College is conducting survey research, designed in collaboration with sociologists, concerning the perceptions of philosophy of students in introductory philosophy courses.}  Class discussion was not in general dominated by confident mathematically-minded men, and this appeared to allow some of the quieter women to speak up and develop their skills and ideas. Some of these women made positive comments to me about feeling more confident or performing better than they had expected---some said this in person, some on their evaluation forms for the course. Overall marks were good and appeared to be slightly better than other instructors' groups, but this difference was not systematically investigated either. It would of course be ideal to see further research into the matter which would produce reliable data. But overall I would recommend all of these strategies to instructors in formal logic courses, since students' responses were positive, their results were good, and my efforts appeared to improve the atmosphere of the exercise classes.

\section{Introducing Women to Role Models: Assigning Original Texts by Female Authors}Another component to the under-representation of women in formal and philosophical logic is that the paucity of women in the field appears to be self-perpetuating. We often hear people who are concerned about the gender gap in logic speculate that female students are not drawn to the field because they lack female role models. I attempted to remedy the lack of female role models while working as a one-on-one supervisor for several colleges (two co-ed, two women-only) on the University of Cambridge second-year Logic course mentioned above. This is primarily a philosophical logic course with a much smaller---and optional---formal component from the one I taught at  St Andrews. I have already pointed out how minuscule the number of women on the official reading list was; in my experience this is not at all unusual, so Cambridge is not noticeably worse in this respect than other universities. The problem of underrepresentation is just very pervasive. 
 My strategy here was assigning some texts from outside the reading lists so that every assignment had at least one piece of set reading by a female author, who whenever possible was presented as the originator of an important view of her own. This strategy was the one that resulted in the most dramatic positive effect that I have seen in my teaching. \emph{All} of my female students, completely unprompted, reacted extremely positively, giving very explicit feedback explaining that it made them feel more included and that it made them enjoy logic more. 

\subsection{My Strategy to Increase the Proportion of Set Reading by Women} The format of the course included weekly short essay-writing assignments on topics such as reference and descriptions, logical form, theories of meaning, interpretations of quantification, theories of truth, intuitionism, and non-classical logics.  I consistently assigned, for each essay, some work by a female logician or philosopher of logic. Works by women comprised at least 25\% of the total reading list for each essay. For a few topics it was not possible to find an appropriately accessible piece of original research by a woman, but for the majority of the assignments at least one female author was presented as the originator of an important view of her own---not as a survey author, but as a towering figure in her own right. Because of the overall paucity of female authors in the field, a typical essay question had assigned reading by two or three male authors and one female author.  Nevertheless this move particularly enthused the female students. As before, to avoid stereotype threat I did not explicitly draw attention to the genders of the students or the authors; I simply informed them of the essay question and the reading list I had set every week. Even before they were asked to give any feedback on the course, female students told me emphatically that they felt more included and that it made them enjoy logic more.

\subsection{Summary of Positive Effects} Female students were very vocal in their appreciation of the original work by female authors, to the extent that I was amazed at the outpourings of enthusiasm. Not only did \emph{all} of the female students, without any kind of prompting from me, comment positively to me about the selection of reading material, but several of them also spoke to their Directors of Studies (Cambridge college subject-specific convenors) who in turn emailed me to pass on even more positive feedback. In previous years I had made extensive use of textbooks and introductory papers by female authors in comparable courses for first- and second-year students, but this never drew a noticeable amount of positive comments from any of the students. Part of the difference may have been that women were consistently included in reading lists for every assignment. But I believe that the main difference here was that female authors were not presented primarily as commentators on the work of (mostly) men, but as serious thinkers in their own right. Many of the positive comments from female students touched upon this aspect. They were very explicit that they appreciated seeing a female author presented as one of the main experts, as someone whose view they were invited to engage with and take seriously. Several also explicitly said that as a result they liked logic a lot more, and felt that logic was `for them' as a result, though they had not previously felt that way.

 \section{Proposals for Gender-Inclusive Tools for Teaching Logic} Given that the most positive results in my efforts to promote gender-inclusive practices in my logic teaching came from including more original texts by female authors among the set reading, my main proposal is to build reading lists, syllabuses, and textbooks around the ideas of women. This proposal divides into various smaller sub-proposals making practical suggestions on how to achieve this.

 \subsection{Building Syllabuses and Reading Lists Around Original Work by Female Authors} The current trend in teaching philosophical logic at the university level, turning away from using textbooks and towards assigning individual papers which students find online, should make it much easier to create more gender-balanced reading lists. This is more difficult in formal logic courses since they rely heavily on textbooks, very few of which are female-authored. Philosophy of logic textbooks are also mostly written by men, and philosophy of logic anthologies are generally hugely male-dominated too. But it is easy for reading lists to be updated to include more journal papers, book chapters and books by female authors. Another component of the strategy, though, must be to make the female authors' work central to the assignment. They must not be featured as an afterthought, nor as handmaidens who only comment on or expound the works of the great men. They should be presented as making interesting original contributions to the debate. 
 
 Some proposals for including more women in syllabuses are very firm that every reading list should be composed of 50\% female authors. I think that, sadly, because logic and philosophy of logic have historically been such male-dominated fields, this is not always feasible in this particular area. But it should be realistic to aim for one-quarter or one-third female authors on the majority of topics in this area. It is also heartening that female students appear to feel much more included even with this much female representation. What is more important than fully proportional representation, it seems, from what I've seen, is that women are presented as serious experts in their fields on a par with the male giants. 
 
 \subsection{Making Essay Questions More Female-Focused}The way essay questions are formulated, or introductory paragraphs on reading lists are written, is also a factor here. The questions and descriptions of the assigned works should present women as originators of serious views as well, and not gravitate only towards the views of the familiar male authors. I frequently used an essay structure where students were asked to compare two opposing views and adjudicate between them, where at least one of these points of view was a female author's.

 Examples include: 
 \begin{itemize}
 \item`Compare and contrast Quine's critique of quantified modal logic with Barcan Marcus' defence of it.' \item`Whose view of conditionals is more successful, Edgington's or Jackson's?'
 \item`Does Grover's pro-sentential theory of truth eliminate the need for correspondence to a fact?'
 \item`Should we be classical or non-classical logicians? In your answer, refer to Dummett and Haack.' 
\end{itemize}

 \subsection{Creating Accessible Archives of Women's Writing in Formal and Philosophical Logic} Another very helpful development would be accessible, searchable archives or databases of papers by women for logic lecturers and tutors to refer to. Many instructors who are sympathetic to the goal of assigning more work by women simply do not know where to find such works, or are unsure where to start looking. This leads to them falling back on the familiar male, pale and stale authors. 
 
 Some efforts have been made to set up databases of writing by women, but generally the resources for them to be properly maintained are not there, and as a result they are not being added to. For instance, http://women.aap.org.au/papers/areas/logic.html only has one paper on logic or philosophy of logic by a woman. A new database of this sort would be very welcome. Alternatively, the editors of existing compendia of philosophical papers could perhaps be persuaded to add a feature that identifies papers by female authors. 
 
 \section{Conclusion} In summary, most of my strategies for addressing hypothesised causes for male dominance in the logic classroom had some positive effect. I attempted to move students' thinking away from the stereotype that men have more of an aptitude for logic or mathematics, but without invoking stereotype threat, by emphasising effort over talent, and praising students for their efforts rather than for innate ability.  There were some positive questionnaires as a result and several students reported that their results were above expectations (anecdotally). But the most noticeable positive effect came from assigning a good proportion of original research by female authors as set reading, with the majority of female students reporting they felt more included as a result. My proposal is therefore to build syllabuses, textbooks and essay questions more specifically around the ideas of women, e.g. Barcan Marcus on reference, modality, and quantification, Blanchette on models and consequence, Edgington on conditionals and knowability, Grover on truth, Haack on non-classical logic, Weiner on Frege.

\subparagraph*{Acknowledgements}

I would like to thank Shannon Dea and the audience at the \emph{Inclusiveness} panel which she organised at the 2015 APA-Central meeting in February 2015. I am grateful to Sophia Connell, too, for her helpful feedback as a Director of Studies of several of the Cambridge colleges where I taught. Thanks are also due to Helen Beebee, Saba Fatima, Kristina Gehrman, and Suzanne Harvey for further discussion. 

\bibliography{bibliography}




\newpage
\thispagestyle{empty}
{\ }

\end{document}